\documentclass[11pt,a4paper,leqno]{article}
\usepackage{a4wide}
\usepackage{latexsym}
\usepackage{amsmath}
\usepackage{amssymb}
\usepackage{graphicx}
\usepackage{xcolor}
\newtheorem{defin}{Definition}
\newtheorem{lemma}{Lemma}
\newtheorem{prop}{Proposition}
\newtheorem{theo}{Theorem}
\newtheorem{corol}{Corollary}
\newtheorem{example}{Example}
\newenvironment{proof}{\medskip\par\noindent{\bf Proof}}{\hfill $\Box$
\medskip\par}

\newcommand{\C}{\mathbb{C}}
\newcommand{\N}{\mathbb{N}}
\newcommand{\Q}{\mathbb{Q}}
\newcommand{\R}{\mathbb{R}}

\begin{document}
\title{On the solutions to linear systems of moment differential equations with variable coefficients}
\author{{\bf A. Lastra}\\
Universidad de Alcal\'{a}, Departamento de F\'{i}sica y Matem\'{a}ticas,\\
Ap. de Correos 20, E-28871 Alcal\'{a} de Henares (Madrid), Spain,\\
{\tt alberto.lastra@uah.es}
}
\date{}
\maketitle
\thispagestyle{empty}
{ \small \begin{center}
{\bf Abstract}
\end{center}

The existence and analyticity of solutions to linear systems of moment differential equations with analytic coefficients is studied. The relation of solutions of such systems with respect to linear moment differential equations is stablished, comparing classical results with the general situation of moment differentiation. 

}

\medskip

\noindent Key words: linear system, moment differential systems, moment differential equations, moment derivative, strongly regular sequence.\\ 2010 MSC: 30D15, 34M03, 34A08.

\bigskip \bigskip

\section{Introduction}

This work contributes to the knowledge of the solutions to linear systems of moment differential equations of the form
\begin{equation}\label{e0}
\partial_{m}y(z)=A(z)y(z)+b(z),
\end{equation}
where $A(z)$ and $b(z)$ are a $n\times n$ matrix and a $n\times 1$ vector of analytic functions on a common domain $D\subseteq \C$, respectively, for some positive integer $n$. $y=y(z)=(y_1(z),\ldots,y_n(z))^{T}$ stands for the vector of unknown functions, and $\partial_my(z)=(\partial_m(y_1),\ldots,\partial_m(y_n))^{T}$, where $m$ is some sequence of positive real numbers, and $\partial_m$ stands for the moment derivative associated with the sequence $m$ (see~(\ref{e144})). This work is a continuation of~\cite{lastra} where the matrix of coefficients is assumed to be constant. 

\vspace{0.3cm}

The notion of moment derivation and the study of moment functional equations was recently firstly developed by W. Balser and M. Yoshino in~\cite{by}, in 2010. This concept generalized that of usual derivative and has a direct reading in terms of $q-$differentiation and Caputo fractional derivatives when chosing $m$ appropriately. This versatility motivates an increasing attention on it in recent years in the scientific community.

The moment derivation, seen as a linear operator, was initially defined on the space of formal power series of coefficients belonging to some Banach space. This operator can be extended to the space of generalized sums and generalized multisums of formal power series (see~\cite{lamisu,lamisu2}). Results on the summability of formal solutions to singularly perturbed partial differential equations in which the coefficients are generalized sums have also been achieved in~\cite{lms}, and on the multisummability of formal solutions in ultraholomorphic classes of functions~\cite{jkls}.

Summability results of formal solutions to moment differential and integral equations have been achieved such as~\cite{michalik12b} and also on Maillet-type results for moment functional equations providing estimates for the coefficients of the formal solutions to such equations~\cite{lamisu3,su} and Stokes phenomenon~\cite{mitk}. 

Regarding the study of systems of moment differential equations, the preceeding work~\cite{lastra} describes all (entire) solutions to homogeneous systems of moment partial differential equations of first order in which the matrix of coefficients is constant, i.e. $A(z)\in\C^{n\times n}$ and $b(z)\equiv 0$ in (\ref{e0}). In that work, the growth at infinity of such solutions is also studied from different points of view. The present work is a natural continuation of~\cite{lastra}.   

The importance in applications of first order linear systems of moment differential equations becomes manifest in different concretizations of the moment sequence. In the case that $m=(\Gamma(1+\alpha p))_{p\ge0}$, the problem can be read in the form of a linear system of fractional differential equations involving Caputo derivatives. The general theory of systems of the form
\begin{equation}\label{e80}
{}^C D^{\alpha}_{z}y(z)=A(z)y(z)+b(z),
\end{equation}
for some $0<\alpha<1$, can be found in detail in the real variable framework in~\cite{brt}, Section 3. In that work, the authors state the existence and unicity, under given initial data, of a system of the form (\ref{e80}), and as a consequence of the non-linear case in~\cite{htrbm}. In this previous work, the existence and uniqueness of solutions to homogeneous linear equations of fractional differential equations of Caputo type of higher order is also analyzed, related to the previous problem. This relationship is also analyzed in the present work. From our point of view, the results obtained while considering general sequences, embracing as a particular case fractional differential equations, despite their versatility,  pay the  price of generality in the form of limited achievements. This will be manifested in the present work. As an example, the absence of a Leibniz-like rule for the moment derivative of a product disables to follow the proofs based on this property, such as the classical variation of constants formula. In practice, the interest of linear, semilinear and linearized systems of fractional differential equations is shown in recent applications such as~\cite{aghror,maab,ucoz} .

In the precedent work~\cite{lastra}, the matrix $A(z)$ in (\ref{e0}) is assumed to be constant and $b\equiv 0$. The main results in that work provide general solution to that problem, which turn out to be entire functions. In the present work, the domain of holomorphy of the solutions to (\ref{e0}) is determined by that of the components of $A(z)$ and $b(z)$. In the study, we consider such domains to be discs centered at the origin. This restriction will be clarified during the work. We have decided to present the results in a progressive manner to emphasize the variation of the radius of convergence of the analytic solution at any point of holomorphy depending on the nature of the sequence of moments considered. Roughly speaking, sequences $m$ whose growth is larger than that of the sequence of factorials (satisfying Assumption (A)) produce no variation on the radius of convergence, whereas sequences $m$ whose growth is that of some root of the sequence of factorials (satisfying Assumption (B)) may produce some variation on the domain of holomorphy. As a matter of fact, most sequences $m$ appearing in applications belong to the set of strongly regular sequences, encompassed in one of the previous situation (Proposition~\ref{prop0}). It is also explained the reason for restricting our study to discs near the origin, avoiding more general results which imply analytic continuation of the solutions to simply connected regions, as it is done in the classical theory of systems of differential equations. As mentioned above, the relationship between first order linear systems and linear equations of higher order is well-known in the classical case, and has also been considered in the framework of fractional differential setting (see~\cite{boritr} and the references therein). The transformation of linear systems of first order moment differential equations into/from linear higher order moment differential equations is considered. In this respect, different situations in which an equation can be transformed into an equivalent system is provided (Theorem~\ref{teo4}). 

The paper is structured as follows. Section~\ref{sec1}  recalls the main facts about strongly regular sequences and generalized summability which will provide the main tools in the study, together with a brief review on fractional calculus. Section~\ref{sec2} is devoted to the study of local existence and uniqueness of Cauchy problems associated to linear systems of the form (\ref{e0}) under different assumptions (Assumption (A) or Assumption (B)) on the sequence $m$. A first result (Theorem~\ref{teo1}) embraces both situations, which hold in the case that $m$ is a strongly regular sequence. The case of nonhomogeneous systems is also considered. In Section~\ref{sec4}, we study the relationship between first order linear systems and linear moment differential equations of higher order. In this respect, the transformation of an equation in the form of a system is guaranteed under certain assumptions (Theorem~\ref{teo4}). The results obtained in~\cite{lastra} are adapted to the case of higher order linear moment differential equations with constant coefficients. 

\vspace{0.3cm}

\textbf{Notation:}

$\N$ denotes the set of nonnegative integers and $\N_0=\N\cup\{0\}$.

$\mathcal{O}(U)$ stands for the set of holomorphic functions on an open set $U\subseteq\C$. $\mathcal{O}\{z\}$ stands for the vector space of holomorphic functions on some neighborhood of the origin.
Given $d\in\R$ and $\theta>0$, we write $S_{d}(\theta)$ for the infinite sector of the Riemann surface of the logarithm $\mathcal{R}$ consisting of all elements $z\in\mathcal{R}$ such that $|z-d|<\theta/2$. We write $D(z_0,r)$ for the open disc centered at $z_0$ and positive radius $r$.

\section{Background on generalized summability and fractional calculus}\label{sec1}

The aim of this section is twofold. On the one hand, we recall the main tools and results appearing in the theory of generalized summability, based on the construction of kernels for generalized summability, departing from strongly regular sequences admitting a nonzero proximate order. The following results and constructions on generalized summability can be found in detail in the seminal work~\cite{sanz}, by J. Sanz. On the other hand, we briefly go through some of the main elements describing fractional calculus that will be appear in the sequel.

\subsection{Strongly regular sequences and generalized summability}

The concept of strongly regular sequence, generalizing Gevrey sequences, is due to V. Thilliez~\cite{thilliez}. It is closely related to function spaces with derivatives bounded in terms of the elements of such sequences.

\begin{defin}
 Given a sequence $\mathbb{M}=(M_p)_{p\ge0}$ of positive real numbers, with $M_0=1$, we say that $\mathbb{M}$ is a strongly regular sequence if the following properties hold:
\begin{itemize}
\item[(lc)] $\mathbb{M}$ is logarithmically convex: $M_p^2\le M_{p-1}M_{p+1}$ for all $p\ge1$.
\item[(mg)] $\mathbb{M}$ is of moderate growth: there exists a positive constant $A_1>0$ such that $M_{p+q}\le A_1^{p+q}M_pM_q$ for all $(p,q)\in\N_0$.
\item[(snq)] $\mathbb{M}$ is non-quasianalytic: there exists $A_2>0$ with $\sum_{q\ge p}\frac{M_q}{(q+1)M_{q+1}}\le A_2\frac{M_p}{M_{p+1}}$, for $p\in\N_0$.
\end{itemize}
\end{defin}

The most outstanding example of a strongly regular sequence is that of Gevrey sequences of a positive order $\alpha>0$, defined by $\mathbb{M}_{\alpha}=(\alpha!^{\alpha})$, of great involvement in the theory of summability of ordinary and partial differential equations, see~\cite{ba2}. Other strongly regular sequences also appear in the study of other functional equations, such as those determining the 1+ level when dealing with difference equations, see~\cite{im1,im2}. The so-called $q-$Gevrey sequence, crucial in the study of $q-$difference equations, satisfies (lc) and (snq), but not (mg) condition, when fixing $q>1$.  

The function $M:[0,\infty)\to[0,\infty)$, defined by $M(0)=0$ and 
\begin{equation}\label{e98}
M(t)=\sup\log_{p\ge0}\left(\frac{t^p}{M_p}\right),\quad t>0
\end{equation}
is well-defined for strongly regular sequences, allowing the construction of kernel functions in the study of generalized summability. As a matter of fact,  the existence of nonzero functions whose asymptotic expansion is zero, is related to the limit opening of a sectorial region of opening $\pi\omega(\mathbb{M})$ where such functions can be defined, where
$$\omega(\mathbb{M})=\frac{1}{\rho(\mathbb{M})}=\left(\lim\sup_{r\to\infty}\max\left\{0,\frac{\log(M(r))}{\log(r)}\right\}\right)^{-1}.$$
For example, it holds that $\omega(\mathbb{M}_{\alpha})=\alpha$.

The main element providing tools for the construction of a theory of generalized summability is the construction of pairs of kernel functions $(e,E)$ whose existence is guaranteed whenever $\mathbb{M}$ admits a nonzero proximate order. 

\begin{defin}[Definition 4,~\cite{jss}]
Let $\mathbb{M}=(M_p)_{p\ge0}$ be a sequence of positive real numbers which satisfies (lc) condition, and $\lim_{p\to\infty}\frac{M_{p+1}}{M_p}=+\infty$. The sequence $\mathbb{M}$ admits the proximate order $\rho:(c,\infty)\to[0,\infty)$ (for some $c\in\R$) if the following statements hold: 
\begin{itemize}
\item[(i)] $\rho$ is a  continuously differentiable function,
\item[(ii)] $\rho(t)\to\rho\in\R$, for some $\rho\in\R_{+}$, 
\item[(iii)] $r\rho'(r)\ln(r)\to 0$,

The previous three conditions describe a proximate order.

\item[(iv)] there exist $A,B>0$ such that $A\le \ln(r)(\rho(r)-d_{\mathbb{M}}(r))\le B$, for large enough $r>0$, where $d_{\mathbb{M}}(r)=\log(M(r))/\log(r)$.
\end{itemize}
\end{defin}

\begin{defin}
Let $\mathbb{M}$ be a strongly regular sequence. We assume that $\omega(\mathbb{M})\in(0,2)$. We say that $(e,E)$ is a pair of kernel functions for $\mathbb{M}$-summability if the following properties hold:
\begin{enumerate}
\item $e\in\mathcal{O}(S_{0}(\omega(\mathbb{M})\pi)$ with $e(z)/z$ being locally uniformly integrable at the origin, i.e. there exists $t_0>0$ and for all $z_0\in S_0(\omega(\mathbb{M})\pi)$ there exists $r_0>0$ (depending on $z_0$) suh that $D(z_0,r_0)\subseteq S_0(\omega(\mathbb{M})\pi)$ with $\int_{0}^{t_0}\sup_{z\in D(z_0,r_0)}|e(t/z)|\frac{dt}{t}\in (0,\infty)$. For every $\delta>0$ there exist $k_1,k_2>0$ such that $|e(z)|\le k_1\exp(-M(|z|/k_2))$, for $z\in S_0(\omega(\mathbb{M})\pi-\delta)$.
\item $e(x)>0$ for $x>0$.
\item $E\in\mathcal{O}(\C)$, with asymptotic growth at infinity determined by
$$|E(z)|\le k_3\exp\left(M(|z|/k_4)\right),\quad z\in\C,$$
for some $k_3,k_4>0$. In addition to this, there exists $\beta>0$ such that for all $\theta\in(0,2\pi-\omega(\mathbb{M})\pi)$ and $R>0$ there exists $k_5>0$ with $|E(z)|\le k_5/|z|^{\beta}$ for $z\in S_{\pi}(\theta)$ with $|z|\ge R$.
\item For every $z\in\C$ one can write
$$E(z)=\sum_{p\ge0}\frac{1}{m(p)}z^{p},$$
with $m=(m(p))_{p\ge0}$ being the sequence of moments (associated to the pair of kernel functions) determined by
$$m(z)=\int_0^{\infty}t^{z-1}e(t)dt,$$
defined on $H=\{z\in\C:\hbox{Re}(z)\ge 0\}$. $m\in\mathcal{O}(\{z\in\C:\hbox{Re}(z)>0\})$ and continuous up to the boundary of $H$.
\end{enumerate}
\end{defin}

Given a strongly regular sequence $\mathbb{M}$ which admits a nonzero proximate order, it is guaranteed the existence of a pair of kernel functions for $\mathbb{M}$-summability (see~\cite{jss,lms}). This is the case of the strongly regular sequences appearing in applications, such as $\mathbb{M}_{\alpha}$ for $\alpha>0$. 

In 2010, W. Balser and M. Yoshino~\cite{by} put forward the definition of a generalized differential operator, acting on formal power series. This notion is associated to a given sequence of positive real numbers $m=(m_p)_{p\ge0}$. The generalized differential operator associated with the sequence $m$ is denoted by $\partial_{m,z}$, or simply $\partial_m$ if the variable can be deduced from the context, and it is known as the moment differential operator due to the sequence $m$ is usually considered to be the moment sequence associated to some positive measure. 

The moment differential operator applied to a formal power series of complex coefficients, $\partial_{m,z}:\C[[z]]\to\C[[z]]$, is defined by
\begin{equation}\label{e144}
\partial_{m,z}\left(\sum_{p\ge0}\frac{a_p}{m_p}z^p\right):=\sum_{p\ge0}\frac{a_{p+1}}{m_p}z^{p}.
\end{equation}
The previous definition can be naturally extended to holomorphic functions on a neighborhood of $z=0$ by identifying the function with its Taylor expansion at the origin. In recent years, the definition of the operator $\partial_{m,z}$ has also been extended to sums and multisums of formal power series in the sense of Poincar\'e, in the framework of summability of formal solutions to functional equations (see~\cite{lamisu,lamisu2}), and also recently applied to solving moment differential equations in the complex domain (see~\cite{lamisu3,su}, among others).

\subsection{Fractional calculus and moment differentiation}\label{fraccal}

This brief section is devoted to review some of the main definitions and results on fractional operators which will be used in the present work. Their importance relies on their relation with moment differentiation. We refer to~\cite{kst} and its references for further details on the topic.

Given $\alpha\in\Q_+$ one can formally define the fractional derivative of order $\alpha$ by
$$\partial_{z}^{\alpha}\left(\sum_{p\ge0}\frac{a_p}{\Gamma(1+\alpha p)}z^{\alpha p}\right)=\sum_{p\ge0}\frac{a_{p+1}}{\Gamma(1+\alpha p)}z^{\alpha p}.$$
The operator $\partial_z^{\alpha}$ is defined on $\C[[z^{\alpha}]]$ with values in the same space of formal power series. For $\alpha=1$ the usual derivative is recovered, whereas given an integer number $k>0$, this operator coincides with Caputo fractional derivative for $1/k$-analytic functions, denoted by ${}^C D^{1/k}_{z}$. It is worth mentioning that 
$$(\partial_{m_{1/k}}\hat{f})(z^{1/k})=\partial_{z}^{1/k}(\hat{f}(z^{1/k})),$$
for every $\hat{f}\in\C[[z]]$, and where $m_{1/k}$ is the moment sequence $(\Gamma(1+p/k))_{p\ge0}$.  The formal inverse of such fractional differential operator coincides with Riemann-Liouville fractional integral, defined by
$$I_{0+}^{\alpha}f(z):=\frac{1}{\Gamma(\alpha)}\int_0^{z}\frac{f(z)}{(z-t)^{1-\alpha}}dt.$$
According to the previous definition, it holds that $I^{\alpha}_{0+}z^{\beta}=\frac{\Gamma(1+\beta)}{\Gamma(+\alpha+\beta)}z^{\alpha+\beta}$. This formal relation can be extended analytically to adequate domains, providing its inverse for sufficiently smooth functions, i.e.:
$$\partial_{z}^{1/k}\circ I_{0+}^{1/k}f=f,$$
in their common domain of definition.

\section{Linear systems of moment differential equations with variable coefficients}\label{sec2}

In this section, we study the existence of solutions of a Cauchy problem determined by a first order system of moment differential equations, with coefficients being holomorphic functions in a disc centered at the origin. Our main result is attained in two steps which allow to illustrate the different nature of their proofs: first, considering a homogeneous system with coefficients defined on a common disc; and second, studying the case of a non homogeneous system. We have decided to maintain all the steps for the sake of readability. 

Let $n\in\N$ with $n\ge 2$ and consider a sequence of positive real numbers $m=(m_p)_{p\ge0}$, with $m_0=1$. Let $A(z)\in\mathcal{O}(D(0,r))^{n\times n}$ be a matrix with components given by holomorphic functions defined in a common disc centered at the origin and positive radius $r$. We consider the system of moment differential equations with initial conditions

\begin{equation}\label{e1}
\left\{ \begin{array}{c}
\partial_my=A(z)y,\\
y(0)=y_0
\end{array}
\right.
\end{equation}
where $y=(y_1,\ldots,y_n)^{t}$ is a vector of unknown functions, with $y_j=y_j(z)$ for $1\le j\le n$, $\partial_my$ denotes the vector $(\partial_my_1,\ldots,\partial_m y_n)^{t}$ and $y_0=(y_{0,1},\ldots,y_{0,n})^{t}\in\C^{n}$.

\vspace{0.3cm}

\textbf{Assumption (A):} There exists $C\in\Q_{+}$ such that $\frac{m_p}{m_{p-1}}\ge \frac{p}{C}$ for every $p\ge 1$.  

\begin{lemma}\label{lema1}
Under Assumption (A) on the sequence $m$, the initial value problem (\ref{e1}) admits a unique solution $y=(y_1,\ldots,y_n)^{t}$, with $y_j\in\mathcal{O}(D(0,r))$ for all $1\le j\le n$.
\end{lemma}
\begin{proof}
Let us write $A(z)=\sum_{p\ge0}A_pz^p\in(\C[[z]])^{n\times n}$. Due to the analyticity property of $A(z)$ it holds that for every fixed $K>1/r$ there exists $c>0$ such that $\left\|A_p\right\|\le c K^p$, where $\left\|\cdot\right\|$ stands for the matrix norm given by the maximum absolute row sum of the matrix. It is always possible to choose $K$ and $c$ to be rational numbers. Let us expand a formal solution of (\ref{e1}) in the form $y(z)=\sum_{p\ge0}y_pz^p$, with $y_p\in\C^n$. We observe that $\partial_m(z^p)=\frac{m_p}{m_{p-1}}z^{p-1}$ for every $p\ge1$ and $\partial_m(1)=0$, and deduce the following recursion formula for the vectors determining the coefficients of the formal solution:
\begin{equation}\label{e81}
y_{p+1}\frac{m_{p+1}}{m_p}=\sum_{k=0}^{p}A_{p-k}y_k,
\end{equation}
for $p\ge0$. Therefore, in the case that the formal solution defines a holomorphic solution of (\ref{e1}), then unicity is guaranteed due to the uniqueness of the coefficients defined from the recursion formula (\ref{e81}), being its first element given by the Cauchy data. We consider the vector norm in $\C^{n}$ defined by $\left\|(z_1,\ldots,z_n)^{t}\right\|=\sup_{1\le \ell\le n}|z_{\ell}|$.


Taking vector norms at both sides of the previous equality and from classical estimates, one concludes that 
$$\left\|y_{p+1}\right\|\frac{m_{p+1}}{m_p}\le \sum_{k=0}^{p}\left\|A_{p-k}\right\|\left\|y_k\right\|\le  c \sum_{k=0}^{p}K^{p-k}\left\|y_k\right\|.$$ We write the previous inequality in the form
\begin{equation}\label{e79}
(p+1)Y_{p+1}\le cK^{p}\left\|y_0\right\|+c\sum_{k=1}^{p}K^{p-k}Y_{k}\frac{m_{k-1}k}{m_k},
\end{equation}
with $Y_{p}=\frac{m_p\left\|y_p\right\|}{m_{p-1}p}$, for every $p\ge 1$, and define $Y_0=\left\|y_0\right\|$.

In this situation, taking into account Assumption (A), it holds that (\ref{e79}) can be estimated by
$$(p+1)Y_{p+1}\le cK^{p}\left\|y_0\right\|+cC\sum_{k=1}^{p}K^{p-k}Y_{k},$$
for every $p\ge1$. We define the sequence $c_0=\left\|y_0\right\|$, and for every $p\ge 1$, 
$$c_{p+1}=\frac{1}{p+1}\left[cK^{p}c_0+cC\sum_{k=1}^{p}K^{p-k}c_{k}\right].$$
It is straight to check that $0\le Y_p\le c_p$ for every $p\ge0$. An analogous argument allows us to construct the sequence $(\tilde{c}_p)_{p\ge0}$ with $\tilde{c}_0=c_0$ and 
\begin{equation}\label{e93}
\tilde{c}_{p+1}=\frac{1}{p+1}\left[c\tilde{C}\sum_{k=0}^{p}K^{p-k}\tilde{c}_{k}\right].
\end{equation}
for every $p\ge1$ and where $\tilde{C}=\max\{1,C\}$, which satisfies that $c_p\le\tilde{c}_{p}$ for every $p\ge0$. We observe the formal power series $\hat{c}(z)=\sum_{p\ge0}\tilde{c}_pz^p$ is a formal solution of the Cauchy problem
\begin{equation}\label{e99}
\left\{ \begin{array}{c}
h'=c\tilde{C}(1-Kz)^{-1}h\\
h(0)=\tilde{c}_0.
\end{array}
\right.
\end{equation}

Therefore, $\hat{c}(z)$ is the Taylor representation of the analytic solution to (\ref{e99}), which is the function $h(z)=c_0(1-Kz)^{-c\tilde{C}/K}$, defining a holomorphic function in $D(0,1/K)$. In view that $Y_p\le c_p\le \tilde{c}_p$ for all $p\ge0$ we conclude that the series $\sum_{p\ge1}\frac{m_p\left\|y_p\right\|}{m_{p-1}p}z^p$ defines a holomorphic function in $D(0,1/K)$. The previous argument is valid for every $K>1/r$. As a conclusion, we obtain that for every $K>1/r$ there exists $c_1>0$ such that 
$$\left\|y_p\right\|\le c_1p\frac{m_{p-1}}{m_p}\frac{1}{K^{p}}\le c_1C\frac{1}{K^{p}},$$
for every $p\ge0$. This entails that $y$ defines a vector of holomorphic functions in $D(0,r)$.
\end{proof} 

\textbf{Remark:} We observe that Assumption (A) is satisfied by every Gevrey sequence $\mathbb{M}_{\alpha}$ of order $\alpha\ge1$, and also by $q-$Gevrey sequences of any positive order, with $|q|>1$.

We state an analogous result under a different assumption.

\vspace{0.3cm}

\textbf{Assumption (B):} There exists a positive constant $C$, and a rational number $0<\alpha<1$ such that $\frac{m_p}{m_{p-1}}\ge \frac{1}{C}\frac{\Gamma(1+\alpha p)}{\Gamma(1+\alpha(p-1))}$ for every $p\ge1$.

\vspace{0.3cm}

Observe that Assumption (A) is the concretion of Assumption (B) for $\alpha=1$. We have decided to split the problem into several parts to remark the viability of the domain of definition of the analytic solution to (\ref{e1}) in these different situations. Assumption (B) involves the use of techniques applied to fractional differential equations which may reduce the radius of convergence of the solution, due to the appearance of Puiseux  series and a natural geometry linked to some multivalued Riemann surface associated to an algebraic function. 

For a given rational number $0<\alpha<1$, we denote by ${}^C D^{\alpha}_{z}$ the Caputo fractional differential operator of order $\alpha$. It is worth recalling that the formal moment derivation $\partial_{m}$ is quite related to the Caputo $\alpha$-fractional differential operator ${}^C D^{\alpha}_{z}$ when fixing the sequence of moments $m=(\Gamma(1+p\alpha))_{p\ge0}$. Indeed, for the previous moment sequence one has that
\begin{equation}\label{e155}
( \partial_{m,z}f)(z^{\alpha})={}^C D^{\alpha}_{z}(f(z^{\alpha}))
\end{equation}
for every $f\in\mathbb{C}[[z]]$. We refer to~\cite{michalik12} (Definition 5 and Remark 1) for further details, and also Section~\ref{fraccal} of the present work.

The proof of Lemma~\ref{lema1} can be adapted to this situation, together with the application of a classical majorant method initialy due to Lindel\"of.

\begin{lemma}\label{lema2}
Under Assumption (B) the initial value problem (\ref{e1}) admits a unique solution $y=(y_1,\ldots,y_n)^{t}$, with $y_j\in\mathcal{O}(D(0,r_0))$ for all $1\le j\le n$, for some $0<r_0\le r$ (depending on $r$, and the elements $C$, $\alpha$ involved in Assumption (B)).
\end{lemma}
\begin{proof}
We sketch the first steps of the proof which follow an analogous scheme as those of Lemma~\ref{lema1}. Let us write $A(z)=\sum_{p\ge0}A_pz^p\in(\C[[z]])^{n\times n}$. For all $K>1/r$ there exists $c>0$ such that $\left\|A_p\right\|\le c K^{p}$. A formal solution $y(z)=\sum_{p\ge0}y_pz^p$ of (\ref{e1}) satisfies the recursion (\ref{e81}), so that the set of inequalities (\ref{e79}) hold for all $p\ge1$. That recursion guarantees unicity of the formal solution as above. From assumption (B) one achieves the inequality 
$$(p+1)Y_{p+1}\le cK^p\left\|y_0\right\|+cC\sum_{k=1}^{p}K^{p-k}Y_{k}\frac{k\Gamma(1+\alpha(k-1))}{\Gamma(1+\alpha k)},$$
for all $p\ge1$. We define the sequence $(c_p)_{p\ge0}$ by $c_0=\left\|y_0\right\|$, and 
$$c_{p+1}=\frac{\Gamma(1+\alpha p)}{\Gamma(1+\alpha(p+1))}\left[cK^pc_0+cC\sum_{k=1}^{p}K^{p-k}c_{k}\right]$$
for all $p\ge1$. It holds that $Y_p p\frac{\Gamma(1+\alpha(p-1))}{\Gamma(1+\alpha p)}\le c_p$ for every $p\ge1$. We define $\tilde{c}_0=c_0$ and $\tilde{c}_{p}$ by (\ref{e93}), for $\tilde{C}=\max\{1,C\}$. It is direct to check that the formal power series $\hat{c}(z)=\sum_{p\ge0}\tilde{c}_pz^p$ is a formal solution of the moment-differential problem 
\begin{equation}\label{e135}
\left\{ \begin{array}{c}
\partial_{m_{\alpha}}h=c\tilde{C}(1-Kz)^{-1}h\\
h(0)=\tilde{c}_0,
\end{array}
\right.
\end{equation}
where $m_{\alpha}=(\Gamma(1+p\alpha))_{p\ge0}$. Regarding the identity (\ref{e155}), the previous problem can be formally rewritten in the form
\begin{equation}\label{e145}
\left\{ \begin{array}{c}
{}^C D^{\alpha}_{z}\tilde{h}=c\tilde{C}(1-Kz^{\alpha})^{-1}\tilde{h}\\
\tilde{h}(0)=\tilde{c}_0,
\end{array}
\right.
\end{equation}
with $\tilde{h}(z)=\hat{c}(z^{\alpha})\in\C[[z^{\alpha}]]$. The recursion formula for the coefficients of the formal Puiseux series solution of (\ref{e145}) coincides with that of (\ref{e93}) and also with the initial data. We now consider the problem (\ref{e145}) in the real variable framework. It is straightforward to check that the formal solution to that problem coincides with that of the integral equation
\begin{equation}\label{e195}
\tilde{h}(x)=\tilde{c}_0+\frac{c\tilde{C}}{\Gamma(\alpha)}\int_0^{x}(x-t)^{\alpha-1}\frac{\tilde{h}(t)}{1-Kt^{\alpha}}dt,
\end{equation}
for $0<x<1/K$. This is a direct consequence of the fact that Caputo derivative is a left inverse of the Riemann-Liouville integral (see Section~\ref{fraccal}). We now consider the sequence of Puiseux polynomials $(\omega_p(x))_{p\ge0}$ defined by $\omega_{p}(x)=\sum_{k=0}^{p}\tilde{c}_{p}x^{p\alpha}$ for every $p\ge0$. We observe that for all $p\ge0$ one has that
\begin{equation}\label{e197}
\omega_{p+1}(x)\le \tilde{c}_0+\frac{c\tilde{C}}{\Gamma(\alpha)}\int_0^{x}(x-t)^{\alpha-1}\frac{\omega_{p}(t)}{1-Kt^{\alpha}}dt.
\end{equation}


Let us fix  
\begin{equation}\label{e201}
0<r_1<\left(\frac{\Gamma(\alpha+1)}{nc\tilde{C}+K\Gamma(\alpha+1)}\right)^{1/\alpha}.
\end{equation}
Observe that $0<r_1<1/K^{1/\alpha}-\epsilon$, for some positive $\epsilon$. We define
$$\Delta=\max\left\{\tilde{c}_0,\frac{\tilde{c}_0}{1-\frac{c\tilde{C}}{\Gamma(\alpha+1)}\frac{r_1^{\alpha}}{1-Kr_1^{\alpha}}}\right\},$$
which is a positive number in view of the choice of $r_1$ in (\ref{e201}). We now prove by an induction argument that 
\begin{equation}\label{e209}
\sup_{x\in[0,r_1]}|\omega_{p}(x)| \le \Delta.
\end{equation}

The estimate (\ref{e209}) is clear for $p=0$. We assume (\ref{e209}) holds up to some positive integer $h$. Regarding (\ref{e197}) one gets for all $0<x<r_1$ that
\begin{align*}
|\omega_{h+1}(x)|& \le \tilde{c}_0+\frac{c\tilde{C}}{\Gamma(\alpha)}\int_0^{x}(x-t)^{\alpha-1}\frac{\omega_{h}(t)}{1-Kt^{\alpha}}dt\\
&\le \tilde{c}_0+\frac{c\tilde{C}}{\Gamma(\alpha)(1-Kr_1^{\alpha})}\left(\sup_{0\le t\le r_1}|\omega_{h}(t)|\right)\int_0^{x}(x-t)^{\alpha-1}dt\\
&\le \tilde{c}_0+\frac{c\tilde{C}}{\Gamma(\alpha+1)(1-Kr_1^{\alpha})}\Delta r_1^{\alpha}\le \Delta,
\end{align*}
regarding the definition of $\Delta$. The series $\omega(x)=\sum_{p\ge0}\tilde{c}_{p}x^{p\alpha}$ converges absolutely on the compact sets of $[0,r_1)$ and is a solution of (\ref{e195}). At this point, one can conclude that for every $0<r_1<\tilde{r}_1$ there exists $C_1>0$ with $\tilde{c}_{p}\le C_1(1/\tilde{r}_1)^{p/\alpha}$. Finally, Assumption (B) guarantees that

$$  \frac{1}{C}\left\|y_p\right\|\le \left\|y_p\right\|\frac{m_p \Gamma(1+\alpha(p-1))}{m_{p-1}\Gamma(1+\alpha p)}=Y_p\frac{p\Gamma(1+\alpha(p-1))}{\Gamma(1+\alpha p)}\le c_p\le \tilde{c}_p\le C_1\left(\frac{1}{\tilde{r}_1}\right)^{p/\alpha}.$$
The conclusion follows from here, with $r_0=\frac{1}{c\tilde{C}/\Gamma(1+\alpha)+1/r}$.
\end{proof}

\vspace{0.3cm}

\textbf{Remark:} We observe that Assumption (B) is satisfied by the sequence $(\Gamma(1+\alpha p))_{p\ge0}$, for $0<\alpha<1$. This corresponds to the classical sequence of moments associated to the  Gevrey sequence $\mathbb{M}_{\alpha}$.

\begin{theo}\label{teo1}
Let $\mathbb{M}=(M_p)_{p\ge0}$ be a strongly regular sequence which admits a pair of kernel functions for $\mathbb{M}$-summability, and let $m=(m(p))_{p\ge0}$ be an associated moment sequence. Then, there exists $R>0$ such that the initial value problem (\ref{e1}) admits a unique solution $y=(y_1,\ldots,y_n)^{t}$, with $y_j\in\mathcal{O}(D(0,R))$ for all $1\le j\le n$. 
\end{theo}
\begin{proof}
We recall from Proposition 5.8 and Remark 6.6~\cite{sanz} that $\mathbb{M}$ and $m$ are both strongly regular sequences which are equivalent sequences in the sense that there exist $c_1,c_2>0$ such that $c_1^pM_p\le m(p)\le c_2^pM_p$, for every $p\ge0$. The fact that $\mathbb{M}$, $m$, $(p!)_{p\ge0}$, and $(\Gamma(1+\alpha p))_{p\ge0}$ for every $\alpha>0$ are strongly regular sequences, and consequently satisfy (lc) and (mg) conditions, one concludes following the argument of Proposition 2.9 (ii)~\cite{sanz}, that:
\begin{itemize}
\item[(i)] in the situation that one can guarantee the existence of $c>0$ such that $m(p)/m(p-1)\ge cp$ for all $p\ge1$, then Assumption (A) holds.
\item[(ii)] in the situation that one can guarantee the existence of $0<\alpha<1$ and $c>0$ such that $m(p)/m(p-1)\ge c\Gamma(1+\alpha p)/\Gamma(1+\alpha (p-1))$ for all $p\ge1$, then Assumption (B) holds.
\end{itemize}
In view of Corollary 1.3(a)~\cite{pet}, there exists $\beta>0$ such that $M_p\ge C_1A_1^pp!^{\beta}$, i.e. there exists $c>0$ such that $M_{p}/M_{p-1}\ge cp!^{\beta}/(p-1)!^{\beta}$ for every $p\ge1$. A direct application of Stirlings formula allows us to conclude that condition (i) or condition (ii) above hold, which allows us to apply Lemma~\ref{lema1} in case that (i) is satisfied, or Lemma~\ref{lema2} in the case that (ii) holds. The result follows directly from here. 

Unicity is guaranteed from the recursion formula satisfied by the coefficients of the formal solution.
\end{proof}

\vspace{0.3cm}

\textbf{Remark:} As a matter of fact, the previous result can be stated in a more general framework in which $m$ is any strongly regular sequence, and not the moment sequence associated to another equivalent strongly regular sequence. We have decided to maintain the previous formulation due to its importance in applications. 

\vspace{0.3cm}

In view of the proof of the previous result, one can deduce the following result.

\begin{prop}\label{prop0}
Let $m$ be a strongly regular sequence. Then, Assumption (A) or Assumption (B) hold.
\end{prop}

\begin{example}
Let $\alpha$ be a positive rational number, and let $r>0$. We consider the equation ($n=1$) 
\begin{equation}\label{e251}
\partial_{m_\alpha}y=\frac{r}{r-z}y,
\end{equation}
under initial condition $y(0)=y_0\in\C$. Here, $m_{\alpha}=(\Gamma(1+\alpha p))_{p\ge0}$. We observe that $A(z)=\frac{1}{r-z}$ is holomorphic on the disc $D(0,r)$. 

First, assume that $\alpha\ge 1$. Stirling's formula determines that
$\frac{\Gamma(1+\alpha p)}{\Gamma(1+\alpha (p-1))}$ has the same asymptotic behavior with respect to $p$ at infinity as $\alpha^{\alpha}(p-1)^{\alpha}$. This entails Assumption (A) holds for that $\alpha$. We observe that for the particular case of $\alpha=1$, equation (\ref{e251}) turns into a first order differential equation with solution given by $y(z)=\frac{y_0}{r-z}$. Therefore, the radius of convergence of the solution to (\ref{e251}) is maintained. For $\alpha>1$, the construction of the proof of Lemma~\ref{lema1} determines that the formal solution to the problem, say $y(z)=\sum_{p\ge0}y_pz^p$, which satisfies that
$$|y_{p+1}|\frac{\Gamma(1+\alpha(p+1))}{\Gamma(1+\alpha p)}\le \sum_{k=0}^{p}K^{p-k}|y_k|,$$
for every $p\ge 0$, and where $K=1/r$. It holds that $|y_p|\le \tilde{c}_{p}\frac{p\Gamma(1+\alpha(p-1))}{\Gamma(1+\alpha p)}$, with $\sum_{p\ge0}\tilde{c}_pz^p$ representing a holomorphic function on the disc $D(0,1/K)$. This concludes holomorphy of $y$ in $D(0,r)$.

In the case that $0<\alpha<1$, then Assumption (B) is satisfied for $C=1$. The majorant series $\hat{c}$ determined in the proof of Lemma~\ref{lema2} is such that $\tilde{h}(x)=\hat{c}(x^{\alpha})$ is a solution of the Caputo fractional differential equation 
$${}^C D^{\alpha}_{x}\tilde{h}=\left(1-\frac{x^{\alpha}}{r}\right)^{-1}\tilde{h},$$
with $\tilde{h}(0)=\tilde{c}_0$. In this situation, one gets 
$$r_0=\frac{1}{\frac{1}{\Gamma(\alpha+1)}+\frac{1}{r}}.$$
\end{example}

Once the radius of the solution of (\ref{e1}) has been treated in the previous results, we depart from a strongly regular sequence to  justify the absence of global solutions defined on more general geometric situations.

The difficulty emerges when defining in a closed and simple manner the moment derivation of a formal power series, as a series in powers of $z-z_0$, for some $z_0\in\C^{\star}$. Indeed, one may write for every $p\ge1$ the power $(z-z_0)^{p}$ in the form of a polynomial in powers of $z$ and apply linearity of $\partial_m$. However, this leads to a much more complicated framework than that of the usual derivative.

Given $z_0\in\C$, it is straight to check that for every formal power series $\hat{f}(z)=\sum_{p\ge0}a_pz^p$, one has that the formal power series in powers of $z-z_0$, $(\hat{f})'(z-z_0)$ and $(\hat{f}(z-z_0))'$, coincide. This is the key point when reducing the study of a system of linear differential equations on some neighborhood of $z_0\in\C$ to the origin. Indeed, given $A(z)\in(\mathcal{O}(D(z_0,r)))^{n\times n}$ and $y_0\in\C$, one has that $y(z)$ is a solution to the Cauchy problem
$$\left\{ \begin{array}{c}
y'=A(z)y,\\
y(z_0)=y_0
\end{array}
\right.$$
defined on $D(z_0,r)$, if and only if $\omega(z)=y(z+z_0)$ is a solution to the Cauchy problem
$$\left\{ \begin{array}{c}
\omega'=B(z)\omega,\\
y(0)=y_0
\end{array}
\right.$$
with $B(z)=A(z+z_0)$. This is no longer the case when considering linear systems of moment differential equations, unless the moment sequence is essentially the sequence of factorials.

\begin{theo}
Let $m$ be a sequence of positive real numbers. Let $z_0\in\C^{\star}$. Then, the following statements are equivalent:
\begin{itemize}
\item
\begin{equation}\label{e415}
(\partial_m\hat{y})(z-z_0)\equiv \partial_{m}(\hat{y}(z-z_0))
\end{equation}
for every $\hat{y}\in\C[[z]]$.
\item $m=(Cp!)_{p\ge0}$, for some $C>0$.
\item $\partial_m$ coincides with usual derivation. 
\end{itemize}
\end{theo}
\begin{proof}
The definition of moment derivation guarantees that the second and third statements are equivalent. The implication that the first statement holds provided that the third is valid is straight from the mentioned classical property of usual derivation. It remains to prove that if the identity (\ref{e415}) holds, then $m$ is essentially Gevrey sequence. Assume that (\ref{e415}) holds for $z^p$ for all integer $p\ge1$. Then, one has that
$$(\partial_{m}z^p)(z-z_0)=\frac{m(p)}{m(p-1)}z^p.$$
On the other hand, 
$$\partial_m((z-z_0)^p)=\partial_m\left[\sum_{j=0}^{p}\binom{p}{j}(-z_0)^{p-j}z^j\right]=\sum_{j=1}^{p}\binom{p}{j}(-z_0)^{p-j}\frac{m(j)}{m(j-1)}z^{j-1}.$$
In order that both polynomials coincide, one should have
$$\frac{m(p)(p-1)!}{m(p-1)(j-1)!(p-j)!}=\frac{p!}{j!(p-j)!}\frac{m(j)}{m(j-1)},$$
for all $j=1,\ldots,p$, or equivalently
$$\frac{m(j)}{m(j-1)j}=\frac{m(p)}{m(p-1)p},$$
for all $1\le j\le p$. The previous reasoning is valid for every positive integer $p$. Therefore, the sequence $(\frac{m(p)}{m(p-1)p})_{p\ge1}$ is constant. There exists $C>0$ such that 
$$\frac{m(p)}{m(p-1)p}=C,\quad p\ge 1.$$
Therefore, $m(p)=Cp!$. Observe that one can define $m(0)=C$. This concludes the proof.
\end{proof}





\begin{example}
Let $k$ be a positive integer. Let $A\in\C^{n\times n}$. We consider the system of moment differential equations (\ref{e1}) with $A(z)=z^kA$ and where $y_0\in\C^{n}$. We assume $m=(m_p)_{p\ge0}$ is a strongly regular sequence.

The formal solution to this initial value problem is given by
\begin{equation}\label{e285}
y(z)=\sum_{p\ge0}\left(\prod_{j=1}^{p}\frac{m_{j(k+1)-1}}{m_{j(k+1)}}\right)A^pz^{(k+1)p}y_0.
\end{equation}
An analogous argument as that of the proof of Theorem~\ref{teo1} guarantees the existence of $\alpha,C>0$, such that 
$$\frac{m_{p+1}}{m_p}\ge C (p+1)^{\alpha},\quad p\ge0.$$
This entails that 
\begin{align*}
\sum_{p\ge0}\left(\prod_{j=1}^{p}\frac{m_{j(k+1)-1}}{m_{j(k+1)}}\right)\left\|A^p\right\||z|^{(k+1)p}\left\|y_0\right\|&\le \sum_{p\ge0}\left(\prod_{j=1}^{p}\frac{1}{C(j(k+1))^{\alpha}}\right)\left\|A\right\|^{p}|z^{k+1}|^{p}\left\|y_0\right\|\\
&=\left\|y_0\right\|\sum_{p\ge0}\frac{1}{p!^{\alpha}}\left(\frac{\left\|A\right\| |z^{k+1}|}{(C(k+1))^{\alpha}}\right)^p,
\end{align*}
which is uniformly convergent in the compact sets of $\C$. This entails that $y(z)$ determines a vector of entire functions. This is coherent with Theorem~\ref{teo1}. 

An upper bound for the growth of the entire solution at infinity can also be determined. Indeed, we observe from Taylor expansion of $y(z)$ at the origin in (\ref{e285}) that $y(z)=\sum_{p\ge0}a_pz^p$, satisfies that
$$\left\|a_p\right\|\le C A^{p}p!^{\alpha/(k+1)}$$
for some $C,A>0$, valid for all $p\ge0$. This is a consequence of Stirling's formula which guarantees the existence of $\tilde{C}_j,\tilde{A}_j>0$, $j=1,2$ such that 
$$\tilde{C}_1\tilde{A}_1^{p}((k+1)p)!^{\frac{\alpha}{k+1}}\le p!^{\alpha}\le \tilde{C}_2\tilde{A}_2^{p}((k+1)p)!^{\frac{\alpha}{k+1}},\qquad p\ge0.$$
Therefore, there exist $c_1,c_2>0$ such that
$$\left\|y(z)\right\|\le c_1\exp\left(c_2|z|^{\frac{k+1}{\alpha}}\right),\qquad z\in\C.$$
This last upper bound is a direct consequence of the application of the following result, which also holds in a more general framework with respect to the families of sequences considered, and the shape of function $M$ under Gevrey settings, which can be found in~\cite{gelfandshilov}. It relates the growth of an entire function at infinity and the growth of the coefficients of its Taylor expansion at the origin (or equivalently any other point). 

\begin{prop}[Proposition 4.5.~\cite{komatsu}]
Let $\mathbb{M}=(M_p)_{p\ge0}$ be a strongly regular sequence, and let $f(z)=\sum_{p\ge0}a_pz^p$ be an entire function. Then, there exists $C_1>0$ such that $|f(z)|\le C_1\exp(M(C_2|z|))$ for $z\in\C$ if and only if there exist $D_1,D_2>0$ such that $|a_p|\le D_1D_2/M_p$, for every $p\ge0$.
\end{prop} 

We observe this is coherent with the classical situation in which $m=(p!)_{p\ge0}$ and $\partial_m$ coincides with the classical derivation and $\alpha=1$. Here, $y(z)=\exp(A\frac{1}{k+1}z^{k+1})y_0$.
\end{example}

\subsection{Nonhomogeneous systems}

In this brief subsection, we provide some information concerning the nonhomogeneous system of moment differential equations
\begin{equation}\label{e326}
\partial_my=A(z)y+b(z),
\end{equation}
where $A\in(\mathcal{O}(D(0,r)))^{n\times n}$ for some $r>0$., and $b=(b_1,\ldots,b_n)^{t}$, with $b_j=b_j(z)\in\mathcal{O}(D(0,r))$ for all $1\le j\le n$.

We omit most of the details from Section~\ref{sec2} that can be mimicked in most of the steps here.

Under Assumption (A), a recursion formula in the form of (\ref{e93}) is obtained
$$\tilde{c}_{p+1}=\frac{1}{p+1}\left[c\tilde{C}\sum_{k=0}^{p}K^{p-k}\tilde{c}_k+\tilde{b}_{p}\right],\qquad p\ge 1,$$
where $\tilde{b}_p=\left\|b_p\right\|$ and $b(z)=\sum_{p\ge0}b_pz^{p}$ for $z\in D(0,r)$. The series $\hat{c}(z)=\sum_{p\ge0}\tilde{c}_pz^p$ formally solves the Cauchy problem
\begin{equation}\label{e99b}
\left\{ \begin{array}{c}
h'=c\tilde{C}(1-Kz)^{-1}h+B(z)\\
h(0)=\tilde{c}_0.
\end{array}
\right.
\end{equation}
with $B(z)=\sum_{p\ge0}\tilde{b}_pz^p$ being a holomorphic function in $D(0,r)$. The variation of constants formula yields the existence of a unique holomorphic solution to (\ref{e99b}) in $D(0,r)$ defined by
$$h(z)=h_1(z)\left[\tilde{c}_0+\int_{z_0}^{z}h_1^{-1}(u)B(u)du\right],$$
for every $z\in D(0,r)$, with $h_1$ being a nonzero solution of the homogeneous equation.

Under Assumption (B), one achieves the recursion formula (\ref{e99b}) and the associated Cauchy problem
\begin{equation}\label{e145b}
\left\{ \begin{array}{c}
{}^C D^{\alpha}_{z}\tilde{h}=c\tilde{C}(1-Kz^{\alpha})^{-1}\tilde{h}+B(z^{\alpha})\\
\tilde{h}(0)=\tilde{c}_0,
\end{array}
\right.
\end{equation}
 with $\tilde{h}(z)=\hat{c}(z^{\alpha})\in\C[[z^{\alpha}]]$. This problem coincides with 
\begin{equation}\label{e195b}
\tilde{h}(x)=\tilde{c}_0+\frac{c\tilde{C}}{\Gamma(\alpha)}\int_0^{x}(x-t)^{\alpha-1}\frac{\tilde{h}(t)}{1-Kt^{\alpha}}dt+\frac{1}{\Gamma(\alpha)}\int_0^{x}(x-t)^{\alpha-1}B(t^{\alpha})dt,
\end{equation}
for $0<x<1/K$. The same choice for $r_1$ as in the expression (\ref{e201}) together with the choice 
$$\Delta=\max\left\{\tilde{c}_0,\frac{\tilde{c}_0+\frac{B_1r_1^{\alpha}}{\Gamma(\alpha)}}{1-\frac{c\tilde{C}}{\Gamma(\alpha+1)}\frac{r_1^{\alpha}}{1-Kr_1^{\alpha}}}\right\},$$
allows to conclude that the problem (\ref{e326}) admits a unique solution with holomorphic components in a disc $D(0,r_0)$, for some $0<r_0<r$, in a similar manner to the proof of Lemma~\ref{lema2}.

\vspace{0.3cm} 

\textbf{Remark:} The case of (\ref{e1}) with $A$ being a constant matrix was considered in~\cite{lastra}. The results obtained in this work remain coherent with those, in which the solutions are vectors of entire functions. 

\section{Linear systems of first order moment differential equations and linear higher order moment differential equations}\label{sec4}

\begin{lemma}\label{lema283}
Let $m$ be a strongly regular sequence. Let $r>0$ and $A(z)\in(\mathcal{O}(D(0,r)))^{n\times n}$. The set of solutions to the system of moment differential equations 
\begin{equation}\label{e2}
\partial_my=A(z)y
\end{equation}
is a subspace of $(\mathcal{O}\{z\})^{n}$ of dimension $n$.
\end{lemma}
\begin{proof}
A proof of the result follows standard arguments as those for a system of first order linear differential equations, i.e. for $m=(p!)_{p\ge0}$ and $\partial_m$ being reduced to the classical derivation, so we only give a sketch of it. 

It is straightforward to check that the set of solutions to (\ref{e2}) is a vector space over $\C$. One can construct a basis of this vector space from $n$ linearly independent solutions to (\ref{e1}) with the initial condition $y_0$ appearing in (\ref{e1}) to be the elements of a basis of $\C^n$, say $e_j$ for $1\le j\le n$. The existence of such solutions is guaranteed by Theorem~\ref{teo1}. The coordinates of any solution of (\ref{e2}), say $y$, with respect to the previous basis coincide with those of $y(0)$ in the basis $\{e_1,\ldots,e_n\}$.
\end{proof}

The previous results motivate the following definition.

\begin{defin}
Let $r>0$ and $A(z)\in(\mathcal{O}(D(0,r)))^{n\times n}$. We also fix a strongly regular sequence $m$. A basis of the vector space of solutions to $\partial_my=A(z)y$ is known as a fundamental solution of the system. A matrix whose columns are given by a fundamental solution is known as a fundamental matrix of solutions.
\end{defin}

\begin{lemma}
Under the assumptions of Lemma~\ref{lema283}, let $Y(z)$ be a fundamental matrix of solutions to (\ref{e2}). Then, $\tilde{Y}(z)$ is another fundamental matrix of solutions to (\ref{e2}) if and only if there exists $C\in\C^{n\times n}$ invertible with $\tilde{Y}(z)=Y(z)C$.
\end{lemma}
\begin{proof}
The only if part is a direct consequence of the columns of $Y$ being basis of the vector space of solutions to (\ref{e2}). $C$ is the change of coordinate matrix. On the other hand, observe that given $C\in\C^{n\times n}$ invertible, then 
$$\partial_m(Y(z)C)=\partial_m(Y(z))C=A(z)Y(z)C=A(z)\tilde{Y}(z).$$
\end{proof}

The usual transformation from linear differential equations to differential systems is available in this framework. In this sense, any sequence of positive numbers $m$, and $a_{j}\in\mathcal{O}(D(0,r))$ for $1\le j\le n$, the homogeneous linear equation
\begin{equation}\label{e294}
\partial_{m}^{n}y_1+a_1(z)\partial_m^{n-1}y_1+\ldots+a_{n-1}(z)\partial_my_1+a_n(z)y_1=0
\end{equation}
is transformed into the system $\partial_my=B(z)y$ with 
\begin{equation}\label{e304}
B(z)= \begin{pmatrix} 
    0 & 1 & 0 & \dots  & 0 & 0\\
	  0 & 0 & 1 & \dots  & 0 & 0\\
    0 & 0 & 0 & \dots  & 0 & 0\\
    \vdots & \vdots & \vdots & \ddots  & \vdots & \vdots\\
    0 & 0 & 0 & \dots  & 0 & 1\\
    -a_n(z) & -a_{n-1}(z) & -a_{n-2}(z) & \dots  & -a_2(z) & -a_1(z)								
		\end{pmatrix}
\end{equation}
and where $\partial_m^{j}$ stands for the iterated operator $\partial_{m}$ acting $j$ times on a scalar function, for $1\le j\le n$. All the previous results can be adapted to linear moment differential equations of the form (\ref{e294}), with solution $y\in\mathcal{O}\{z\}$. The Cauchy data associated to such a problem are $y(0)=y_0$, $\partial_my(0)=y_0^1$, $\ldots$, $\partial_{m}^{n-1}y(0)=y_0^{n-1}$, for some $y_0^{j}\in\C$ for $0\le j\le n-1$. The next result is a direct consequence of the previous identification.

\begin{corol}
Let $r>0$ and let $m$ be a strongly regular sequence. We take $a_j\in\mathcal{O}\{z\}$ for $1\le j\le n$. The set of solutions to (\ref{e294}) is a subspace of dimension $n$ of the vector space $\mathcal{O}\{z\}$.
\end{corol}

Consequently, one can state the concept of fundamental solution in the framework of linear systems of a positive order.

\begin{defin}
Let $r>0$ and let $a_j\in\mathcal{O}\{z\}$ for $1\le j\le n$. We also fix a strongly regular sequence $m$. A basis of the vector space of solutions to (\ref{e294}) is known as a fundamental solution of the linear moment differential equation. 
\end{defin}

The construction of a fundamental matrix of solutions associated to (\ref{e2}) can be performed in an analogous manner as in the classical case. The determinant of such fundamental matrix satisfies classical non nullity properties which are not described in this work, as they can be deduced in the same way as in the classical theory from appropriate fixed assumptions on the sequence $m$. 

The following result is an example of the previous assertion. Observe that every strongly regular sequence $m=(m_{p})_{p\ge0}$ satisfies the so-called derivability condition $\sup_{p\ge0}\left(\frac{m_{p+1}}{m_p}\right)^{p}<\infty$, present in a wide number of results in the literature. In this respect, see for example~\cite{thilliez0}, and the references therein.

\begin{lemma}\label{lema4}
Let $m$ be a strongly regular sequence. Let $Y(z)\in(\mathcal{O}\{z\})^{n\times n}$. Then, there exists a unique $A(z)\in(\mathcal{O}\{z\})^{n\times n}$ such that $Y$ is the fundamental matrix of the problem (\ref{e2}). 
\end{lemma}
\begin{proof}
The result follows from the classical construction $A(z)=(\partial_m Y(z))Y^{-1}(z)$.
\end{proof}

On the other hand, the classical procedure to transform the system (\ref{e2}) into a linear equation of the form (\ref{e294}) based on the construction of a basis determined by a cyclic vector is no longer available in general, due to the absence of a Leibniz rule with respect to moment derivatives. Some results can be achieved under appropriate assumptions.

\begin{defin}
Let $A\in\C^{n\times n}$. We say $v\in\C^n$ is a cyclic vector associated to $A$ if the set $\{v,vA,\ldots,vA^{n-1}\}$ is a basis of $\C^n$.
\end{defin}

Observe that the existence of a cyclic vector associated to a matrix $A$ is not always guaranteed. For example, for $n\ge 3$ and $\hbox{rank}(A)\le n-2$, then a cyclic vector associated to $A$ can not exist. Indeed, it is a classical result from linear algebra that a cyclic vector associated to $A$ exists if and only if the characteristic and the minimal polynomial of $A$ coincide (up to a constant factor). We recall that the characteristic polynomial of $A$ can be defined by $c_A(x)=\det(A-xI)$ whereas the minimal polynomial of $A$ is a polynomial $c_m$ of minimum degree such that $c_m(A)=0$. It holds that $m_A|c_A$. With respect to Jordan representation of a matrix, the existence of a cyclic vector is equivalent to the fact that each eigenvalue of $A$ can only appear in one Jordan block.
 This holds in particular if $A$ has different simple eigenvalues.



A more practical way to determine that minimal and characteristic polynomials coincide is the following well-known result.

\begin{lemma}
If gcd of all nonzero minors of order $n-1$ of the matrix $A-xI$ is a constant (is a polynomial of degree 0 on $x$), then $c_A=m_A$.
\end{lemma}

\begin{lemma}\label{lema327}
Let $A\in\C^{n\times n}$ be a constant matrix which admits a cyclic vector, and let $m$ be a strongly regular sequence. Then, there exist complex numbers $a_1,\ldots,a_n\in\C$ and an invertible matrix $T\in\C^{n\times n}$ such that if $y(z)$ is solution of the system $\partial_m y=Ay$, and we define $\tilde{y}=Ty$, then $\tilde{y}$ is a solution of $\partial_m y=By$, with $B$ as in (\ref{e304}).
\end{lemma}
\begin{proof}
The classical construction of a basis of $\C^n$ given by $\{h_1,h_1A,\ldots,h_1 A^{n-1}\}$, for some cyclic vector $h_1\in\C^n$ associated to $A$ allows to conclude the result. Observe that the coefficients of $h_1A^n$ in the previous basis determine the opposite of the numbers $a_1$ up to $a_n$.
\end{proof}

\vspace{0.3cm}

\textbf{Remark:} Observe that $y_j\in(\mathcal{O}(\C))^{n}$ and $\tilde{y}_j\in\mathcal{O}(\C)$ for every $1\le j\le n$ in Lemma~\ref{lema327}, in view of Lemma 2~\cite{lastra}.


\vspace{0.3cm}

The following result allows to transform a linear system into a higher order linear moment differential equation for a non-constant matrix $A$ under certain assumptions.

\begin{theo}\label{teo4}
Let $A(z)\in(\mathcal{O}\{z\})^{n\times n}$. Let $m$ be a strongly regular sequence and consider the system of moment differential equations $\partial_my=A(z)y$. Assume that:
\begin{itemize}
\item[(i)] The matrix $A_0:=A(0)$ admits a cyclic vector $v_0\in\C^n$.
\item[(ii)] For every $p\ge1$ and $0\le j\le n-2$, it holds that $v_0A(0)^{j}A^{(p)}(0)\equiv 0$ ($A^{(p)}$ stands for the matrix whose components are the $p-$th derivatives of the components of $A$).
\end{itemize}

Then, there exist $a_1,\ldots,a_n\in\mathcal{O}\{z\}$ and an invertible matrix $T\in\C^{n\times n}$ such that if $y(z)$ is solution of the system $\partial_m y=Ay$, and we define $\tilde{y}=Ty$, then $\tilde{y}$ is a solution of $\partial_m y=B(z)y$, with $B(z)$ as in (\ref{e304}). 

In addition to this, if $A(z)\in(\C_p[z])^{n\times n}$ ($\C_p[z]$ stands for the vector space of polynomials of complex coefficients and degree at most $m$) then $a_j(x)\in\C_p[x]$ for every $1\le j\le n$.
\end{theo}
\begin{proof}
Write $A(z)=\sum_{p\ge0}A_pz^p$ for $z\in D(0,r_0)$, for some $r_0>0$. Let $v_0\in\C^n$ be a cyclic vector of $A_0$. We write $v_1=v_0A(z)=v_0A_0$ in view of (ii). Recursively, we define $v_j=v_{j-1}A(z)=v_0A_0^{j}$ for every $1\le j\le n-1$. Let $T$ be the matrix with rows given by $v_j$ for $0\le j\le n-1$. Observe that $T$ is an invertible constant matrix. We define $B(z)$ as in (\ref{e304}), with its last row given by $v_0A_0^{n-1}A(z)T^{-1}$. Observe that 
$$\partial_{m}T+TA(z)=TA(z)=B(z)T(z)$$
to conclude the proof.

In the case that $A(z)\in (\C_p[z])^{n\times n}$, then the previous construction yields to the polynomial nature of the coefficients in (\ref{e304}).
\end{proof}

\vspace{0.3cm}

\textbf{Remark:} We observe that a more general configuration of the matrix $A_0$ in (i) of Theorem~\ref{teo4} would provide a weaker result in this direction. More precisely, assume that $A_0$ can be written as a direct sum of $m \ge2$ matrices representing invariant subspaces. Then, the system (\ref{e2}) can be rewritten in the form of $m$ linear moment differential equations via a blocked matrix with blocks of the form (\ref{e304}) together with linear first order moment differential equations in the case of invariant subspaces corresponding to diagonalizable terms in the decomposition of $A_0$.  

\vspace{0.3cm}

\begin{example}
We consider the system of moment differential equations $\partial_my=A(z)y$, with
$$A(z):=\begin{pmatrix}
3z+9 & 9z+8 & 6z-11\\
z-3 & z+2 & 2z+7\\
-5z-3 & -5z+8 & -10z+1
\end{pmatrix}.$$
We observe that $v_0=(1,2,1)$ is a cyclic vector of 
$$A_0:=\begin{pmatrix}
9 & 8 & -11\\
-3 & 2 & 7\\
-3 & 8 & 1
\end{pmatrix}$$
due to the set 
$$\{v_0,v_0A_0,v_0A_0^2\}=\{(1,2,1),(0,20,4),(-72,72,144)\}$$
determines a basis of $\C^3$. Therefore condition (i) in Theorem~\ref{teo4} holds. In addition to this, we have that $v_0A_1=v_0A_0A_1=(0,0,0)$ and $v_0A_0^jA^{(p)}(0)=(0,0,0)$ due to $A(z)$ is a matrix polynomial of degree 1, for every $j\ge2$. This entails condition (ii) in Theorem~\ref{teo4}. Following the construction in the proof of Theorem~\ref{teo4} we derive that 
\begin{equation}\label{e392}
v_0A_0^2A(z)T^{-1}=(-1296 z-432,108 z+36,-6 z+12),
\end{equation} 
concluding that the solution to the system of moment differential equations is related to that of $\partial_m\tilde{y}=B(z)\tilde{y}$, where $B(z)$ is as in (\ref{e304}), with its last row given by given by (\ref{e392}), i.e. related to the linear moment differential equation 
$$\partial_{m}^{3}y+(-1296 z-432)\partial_{m}^{2}y+(108 z+36)\partial_{m}y+(-6 z+12)y=0.$$
\end{example}

\vspace{0.3cm}

\textbf{Remark:} A variation of constants formula for a closed expression of the solutions of an nonhomogeneous system 
\begin{equation}\label{e470}
\partial_my=A(z)y+b(z),
\end{equation}
for some $b(z)\in(\mathcal{O}\{z\})^{n}$, seems unlikely to exist due to the absence of a general Leibniz rule considering moment derivation. This will be the object of a future study. 

\vspace{0.3cm}

The results previously obtained in~\cite{lastra} on the asymptotic behavior of the solutions to linear systems of moment differential equations can be naturally adapted to linear moment differential equations of constant coefficients of the form (\ref{e294}). We only sketch the construction and summarize the main results obtained therein. 

\begin{defin}[Definition 5,~\cite{lastra}]
Let $\lambda\in\C$ and $h\ge0$ be an integer. We also fix a sequence of positive real numbers $(m_p)_{p\ge0}$. We define the formal power series
$$\Delta_{h}E(\lambda,z):=\sum_{p\ge h}\binom{p}{h}\frac{\lambda^{p-h}z^p}{m_p}.$$
\end{defin}

\begin{prop}\label{prop1}
Let $\mathbb{M}$ be a strongly regular sequence which admits a nonzero proximate order $\rho(t)\to\rho>0$, and let $E(\cdot)$ be a kernel function for generalized summability associated to it. Then, for every nonnegative integer $h$ and $\lambda\in\C$ the formal power series $\Delta_hE(\lambda,z)$ defines an entire funtion. In addition to this, the following statements hold on its asymptotic behavior at infinity:
\begin{itemize}
\item there exist $C_1,C_2>0$ such that
$$|\Delta_{h}E(\lambda,z)|\le C_1\exp(M(C_2|z|),\qquad z\in\C,$$
where $M(\cdot)$ stands for the function defined in (\ref{e98}).
\item $\Delta_hE(\lambda,z)$ is an entire function of order $\rho$ and type $|\lambda|^{\rho}$, i.e.
$$\rho=\lim\sup_{r\to\infty}\frac{\ln^+(\ln^+(M_{h,\lambda}(r)))}{\ln(r)},$$
$$|\lambda|^{\rho}=\lim\sup_{r\to\infty}\frac{\ln^+(M_{h,\lambda}(r))}{r^{\rho}},$$
with $\ln^+(x)=\max\{0,\ln(x)\}$ for $x>0$, and $M_{h,\lambda}(r):=\max\{|\Delta_hE(\lambda,z)|:|z|=r\}$. 
\end{itemize}
\end{prop}

The results obtained in~\cite{lastra} can be rewritten as follows according to the framework of Cauchy problems involving linear moment differential equations with constant coefficients.

Let $a_j\in\C$ for $1\le j\le n$. We also fix $y_0^{j}\in\C$ for $1\le j\le n-1$. Let $\mathbb{M}$ be a strongly regular sequence admitting a nonzero proximate order $\rho(t)\to\rho>0$, for $t\to\infty$. We consider a pair of kernel functions for generalized summability $(e,E)$ and the corresponding sequence of moments $m$, and define the problem
\begin{equation}\label{e294b}
\left\{ \begin{array}{l}
\partial_{m}^{n}y+a_1\partial_m^{n-1}y+\ldots+a_{n-1}\partial_my+a_ny=0\\
\partial_m^{j}y(0)=y_0^{j}.
\end{array}
\right.
\end{equation}  

We first describe the solution space in a particular case, which puts into light the structure of the solutions of the equation in (\ref{e294b}) as a direct sum of subspaces.

\begin{lemma}\label{lema565}
In the previous situation, we fix a positive integer $\ell$ and $\lambda\in\C$, and consider the moment differential equation
\begin{equation}\label{e567}
(\partial_m-\lambda)^{\ell}y=0.
\end{equation}
Then, the set $\{\Delta_0E(\lambda,z),\Delta_1E(\lambda,z),\ldots,\Delta_{\ell-1}E(\lambda,z)\}$ is a basis of the space of solutions to (\ref{e567}).
\end{lemma}
\begin{proof}
The reformulation of the problem in the form of a first order linear system of moment differential equations guarantees that the solution to (\ref{e567}) is a vector space of dimension $\ell$. Taking into account that $\Delta_0 E(\lambda,z)=E(\lambda z)$ for every $\lambda\in\C$ together with the property that for all $h\ge 1$ and $\lambda\in\C$ one has that
$$(\partial_m-\lambda)(\Delta_hE(\lambda,z))=\Delta_{h-1}E(\lambda,z),\quad z\in\C$$
we get that for all $0\le k\le \ell-1$:
\begin{equation}\label{e571}
(\partial_m-\lambda)^{\ell}(\Delta_kE(\lambda,z))\equiv 0.
\end{equation}
This guarantees that $\Delta_kE(\lambda,z)$ for $0\le k\le\ell-1$ is a solution of (\ref{e567}). We now prove their linear independence. Assume a linear combination 
\begin{equation}\label{e675}
\sum_{k=1}^{\ell}C_k  \Delta_{k-1}E(\lambda,z)\equiv0
\end{equation}
for some constants $C_{k}$, $1\le k\le\ell$. The evaluation at $z=0$ of the previous expression yields $C_1=0$. Taking moment derivatives at (\ref{e675}) one arrives at 
$$0=\sum_{k=2}^{\ell}C_k \partial_m(\Delta_{k}E(\lambda,z))=C_2E(\lambda z)+C_2\lambda \Delta_{1}E(\lambda,z)+\sum_{k=3}^{\ell}C_k (\lambda\Delta_{k-1}E(\lambda,z)+ \Delta_{k-2}E(\lambda,z)).$$
The evaluation of the previous expression at $z=0$  determines $C_2=0$. A recursive argument consisting on taking moment derivatives at the resulting expression allows to conclude, having at every step that 
$$0=C_hE(\lambda z)+zF(\{\Delta_{j}E(\lambda ,z)\}_{1\le j\le \ell}),$$
with $F$ being a polynomial. This concludes the proof.
\end{proof}

\begin{theo}\label{teo6}
In the previous situation, let the characteristic polynomial associated to (\ref{e294b}) be written in the form $(\lambda-\lambda_1)^{\ell_1}\cdots(\lambda-\lambda_{r})^{\ell_r}$ for some $1\le r\le n$, and with $\ell_1+\ldots+\ell_r=n$, $\lambda_j\in\C$ for $1\le j\le r$ with $\lambda_{j_1}\neq\lambda_{j_2}$ if $j_1\neq j_2$. Then, there exist $c_{j,k}$ for $1\le j\le n$ and $1\le k\le \ell_{j}$ such that
$$y(z)=\sum_{j=1}^{r}\sum_{k=1}^{\ell_j}c_{j,k}\Delta_{k-1}E(\lambda_j, z).$$
\end{theo}
\begin{proof}
An analogous argument as that of the beginning of the proof of Lemma~\ref{lema565} allows us to affirm that the set of solutions to the equation in (\ref{e294b}) determines a vector space of dimension $n$. Once provided a basis of solutions, the coordinates in that basis of the unique solution to (\ref{e294b}) are determined by the Cauchy data. This again can be directly observed from the matrix associated to the problem.

We conclude the result by proving the linear independence of the set $\{\Delta_{k-1}E(\lambda_j,z):1\le j\le r,1\le k\le \ell_j\}$. Observe that $\Delta_{k-1}E(\lambda_j z)$ has a zero of order $k-1$ at $z=0$. Given a linear combination 
$$C_{11}\Delta_0E(\lambda_1,z)+\ldots+C_{1\ell_1}\Delta_{\ell_1-1}E(\lambda_1,z)+\ldots+
C_{r1}\Delta_0E(\lambda_r,z)+\ldots+C_{r\ell_r}\Delta_{\ell_r-1}E(\lambda_r,z)=0,$$
we write it in the form
\begin{equation}\label{e698}
\Phi_0(z)+z\Phi_1(z)+\ldots+z^{\xi_r-1}\Phi_{\xi_r-1}=0,
\end{equation}
with $\xi_r=\max_{j=1,\ldots,r}\ell_j$ and where 
$$\Phi_{k-1}(z)=\frac{1}{z^{k}}\sum_{j=1}^{r}C_{jk}\Delta_{k-1}E(\lambda_j,z),$$
for every $1\le k\le \xi_r-1$, defining $C_{jk}=0$ for $k\ge \ell_j$. By applying the classical derivation recurrently on the expression (\ref{e698}), it is straight to check that $\Phi_{k-1}\equiv 0$ for $1\le k\le \xi_r-1$. It only rests to check that for every $k$, the linear combination determining $z^{k-1}\Phi_{k-1}$ is only possible when the coefficients are zero. Such linear combination is given by
$$C_{1k}\Delta_{k-1}E(\lambda_1,z)+ \ldots+C_{rk}\Delta_{k-1}E(\lambda_r,z)=0.$$
In order that the previous equality holds, the first $r$ coefficients in the previous power series vanish. In other words, one has the homogeneous linear system of $r$ equations
$$\left\{\sum_{j=1}^{r}\binom{k-1+h}{h}\frac{1}{m(k-1+h)}\lambda_j^{h}C_{jk}=0\right\}_{h=0\ldots r-1},$$
with associated Vandermonde matrix of coefficients $(\lambda_{j}^{h})_{0\le j,h\le r-1}$, whose determinant is different from zero, due to $\lambda_{j_1}\neq\lambda_{j_2}$ for $j_1\neq j_2$. Therefore, the only solution is given by $C_{jk}=0$ for $1\le j\le r$. The linear independence and therefore the conclusion of the result follow from here.
\end{proof}

As a consequence of Proposition~\ref{prop1}, the following results on the asymptotic behavior on the solution to (\ref{e294b}) are achieved. We omit their proofs which are a direct application of Proposition~\ref{prop1} to Theorem 5, Proposition 4 and Theorem 7 of~\cite{lastra}, respectively. 

\begin{corol}
Under the assumptions of Theorem~\ref{teo6}, it holds that:
\begin{enumerate}
\item if $\lambda=0$ is the only root of the characteristic polynomial associated with (\ref{e294b}), then the solution of (\ref{e294b}) is a polynomial.
\item Let $\{\lambda_j:1\le j\le r\}$ be the roots of the characteristic polynomial associated with (\ref{e294b}). Then, the solution of (\ref{e294b}) is an entire function of order $\rho$ and type upper bounded by $\sigma:=(\max_{1\le j\le r}|\lambda_j|)^{\rho}$.
\end{enumerate}
\end{corol}



\begin{corol}
Under the assumptions of Theorem~\ref{teo6}, suppose that the characteristic polynomial associated to the moment differential equation in (\ref{e294b}) has $n$ simple roots, $\lambda_j\in\C$ for $1\le j\le n$. Then, the unique entire solution to (\ref{e294b}), say $y=y(z)$, satisfies that 
$$|y(re^{i\theta})|\le \frac{C}{r^{\beta}},\qquad r\ge R_0,$$
for some $C,\beta,R_0>0$ provided that $\theta$ belongs to the set
$$\bigcap_{j=1}^{k}\left\{\theta\in\R:\frac{\omega(\mathbb{M})\pi}{2}<\theta+\hbox{arg}(\lambda_j)<2\pi-\frac{\omega(\mathbb{M})\pi}{2}\right\},$$
whenever the previous set is not empty.
\end{corol}

\begin{defin}
Let $\rho(t)$ be a proximate order associated to $f\in\mathcal{O}(\C)$, meaning that $\rho(t)$ is a proximate order and $\sigma_f$ is a positive number, where
$$\sigma_f=\lim\sup_{r\to\infty}\frac{\ln(M_f(r))}{r^{\rho(r)}}.$$
The generalized indicator of $f$ is defined by
$$h_f(\theta)=\lim\sup_{r\to\infty}\frac{\ln|f(re^{i\theta})|}{r^{\rho(r)}}.$$
\end{defin}

\begin{corol}
Under the assumptions of Theorem~\ref{teo6}, suppose that the characteristic polynomial associated to the moment differential equation in (\ref{e294b}) has $n$ simple roots, $\lambda_j\in\C$, for $1\le j\le n$. Then, the solution of (\ref{e294b}), say $y\in\mathcal{O}(\C)$, satisfies that
$$h_{y}(\theta)\le\max\{|\lambda_j|^{\rho}h_{E}(\theta+\arg(\lambda_j)):1\le j\le n\},$$
for all $\theta\in\R$.
\end{corol}

\vspace{0.3cm}

\textbf{Acknowledgements:} A. Lastra is partially supported by the project PID2022-139631NB-I00 of Ministerio de Ciencia e Innovaci\'on, Spain.


\begin{thebibliography}{99}
\bibitem{aghror} R. Agarwal, S. Hristova, D O'Regan, \emph{Mittag-Leffler stability for impulsive Caputo fractional differential equations}, Differ. Equ. Dyn. Syst. 29, No. 3 (2021) 689--705. 
\bibitem{ba2} W. Balser, \emph{Formal power series and linear systems of meromorphic ordinary differential equations.} Universitext. Springer-Verlag, New York, 2000. 
\bibitem{by} W. Balser, M. Yoshino, \emph{Gevrey order of formal power series solutions of inhomogeneous partial differential equations with constant coefficients,} Funkcial. Ekvac. 53 (2010) 411--434. 
\bibitem{boritr} B. Bonilla, M. Rivero, J. J. Trujillo, \emph{Linear differential equations of fractional order}, J. Sabatier (ed.) et al., Advances in fractional calculus. Theoretical developments and applications in physics and engineering. Dordrecht: Springer (2007) 77--91. 
\bibitem{brt} B. Bonilla, M. Rivero, J. J. Trujillo, \emph{On systems of linear fractional differential equations with constant coefficients,} Appl. Math. Comput. 187 (2007) No. 1, 68--78.
\bibitem{caputo} M. Caputo, \emph{Lineal model of dissipation whose Q is almost frequency independent II}, Geophys. J. R. Astronom. Soc. 13 (1967) 529--539.

\bibitem{gelfandshilov} I. M. Gelfand, G. E. Shilov, \emph{Generalized functions}, Vol. 2. Academic New York, 1965. 
\bibitem{htrbm} N. Hayek, J. Trujillo, M. Rivero, B. Bonilla, J. C. Moreno, \emph{An extension of Picard-Lindel{\"o}ff theorem to fractional differential equations}, Appl. Anal. 70, No. 3-4 (1999) 347--361. 
\bibitem{im1} G. Immink, \emph{Accelero-summation of the formal solutions of nonlinear difference equations,} Ann. Inst. Fourier (Grenoble) 61, no. 1  (2011) 1--51.
\bibitem{im2} G. Immink, \emph{Exact asymptotics of nonlinear difference equations with levels 1 and 1+}, Ann. Fac. Sci. Toulouse Math. (6) 17, no. 2, (2008) 309--356.
\bibitem{jkls} J. Jim\'enez-Garrido, S. Kamimoto, A. Lastra, J. Sanz, \emph{Multisummability in Carleman ultraholomorphic classes by means of nonzero proximate orders}, J. Math. Anal. Appl. 472 (1) (2019) 627--686.
\bibitem{jss} J. Jim\'enez-Garrido, J. Sanz, G. Schindl, \emph{Log-convex sequences and nonzero proximate orders}, J. Math. Anal. Appl. 448(2) (2017) 1572--1599.
\bibitem{kst} A. A. Kilbas, H. M. Srivastava, J. J. Trujillo, Theory and applications of fractional differential equations. North-Holland Math. Stud. 204, Elsevier, Amsterdam, 2006.
\bibitem{komatsu} H. Komatsu, \emph{Ultradistributions. I: Structure theorems and a characterization}, J. Fac. Sci. Univ. Tokyo, Sect. I A 20 (1973) 25--105. 
\bibitem{lastra} A. Lastra, \emph{Entire solutions of linear systems of moment differential equations and related asymptotic growth at infinity,} Differ. Equ. Dyn. Syst. (2022).
\bibitem{lms} A. Lastra, S. Malek, J. Sanz, \emph{Summability in general Carleman ultraholomorphic classes}, J. Math. Anal. Appl. 430 (2015) 1175--1206. 
\bibitem{lamisu} A. Lastra, S. Michalik, M. Suwi{\'n}ska, \emph{Summability of formal solutions for a family of generalized moment integro-differential equations,} Fract. Calc. Appl. Anal. 24 (2021) No. 5, 1445--1476. 
\bibitem{lamisu2} A. Lastra, S. Michalik, M. Suwi{\'n}ska, \emph{Multisummability of formal solutions for a family of generalized singularly perturbed moment differential equations}, Results Math. 78, 49 (2023).
\bibitem{lamisu3} A. Lastra, S. Michalik, M. Suwi{\'n}ska, \emph{Estimates of formal solutions for some generalized moment partial differential equations,} J. Math. Anal. Appl. 500(1), 18 (2021)


\bibitem{maab} M. M. Matar, E. S. A. Skhail, \emph{On stability analysis of semi-linear fractional differential systems,} Math. Methods Appl. Sci. 43, No. 5 (2020) 2528--2537. 
\bibitem{michalik12} S. Michalik, \emph{Multisummability of formal solutions of inhomogeneous linear partial differential equations with constant coefficients}, J. Dyn. Control Syst. 18 (2012) 103--133.
\bibitem{michalik12b} S. Michalik, \emph{Analytic solutions of moment partial differential equations with constant coefficients}, Funkcial. Ekvac. 56 (1) (2013) 19--50.
\bibitem{mitk} S. Michalik, B. Tkacz, \emph{The Stokes phenomenon for some moment partial differential euqations}, J. Dyn. Control Syst. 25 (4) (2019) 573--598.
\bibitem{pet} H.-J. Petzsche, \emph{On E. Borel's theorem}, Math. Ann. 282 (1988) 299--313.
\bibitem{sanz} J. Sanz, \emph{Flat functions in Carleman ultraholomorphic classes via proximate orders}, J. Math. Anal. Appl. 415(2) (2014) 623--643.
\bibitem{su} M. Suwi{\'n}ska, \emph{Gevrey estimates of formal solutions for certain moment partial differential equations with variable coefficients,} J. Dyn. Control Syst. 27(2) (2021) 355--370. 
\bibitem{thilliez} V. Thilliez, \emph{Division by flat ultradifferentiable functions and sectorial extensions}, Result. Math. 44 (2003) 169--188. 
\bibitem{thilliez0} V. Thilliez, \emph{Smooth solutions of quasianalytic or ultraholomorphic equations}, Monatsh. Math. 160 (2010) 443--453.
\bibitem{ucoz} E. U\c{c}ar, N. \"Ozdemir, \emph{A fractional model of cancer-immune system with Caputo and Caputo–Fabrizio derivatives,} Eur. Phys. J. Plus 136 (2021) 43.
\end{thebibliography}
\end{document}